\documentclass{birkjour}
\usepackage{amsmath}
\usepackage{amssymb,amsthm}
\usepackage{graphicx}
\usepackage{enumerate}
\usepackage{amsfonts}
\usepackage{pstricks}
\usepackage{verbatim}
\usepackage{bbm}
\usepackage{bbold}
\usepackage{mathrsfs, hyphenat}
 \newtheorem{thm}{Theorem}[section]

 \newtheorem{prop}[thm]{Proposition}
 \theoremstyle{definition}
 \newtheorem{defn}[thm]{Definition}
 \theoremstyle{remark}

 \numberwithin{equation}{section}

\def\int{\mathbb Z}

\def\ZN{{\mathbb Z}_N}
\def\H{{\mathbb H}}
\def\Z2{{\mathbb Z}_2}

\def\E_B{{\frak E}_{{\frak B}}}%
\def\Lin_F{\mbox{Lin}(\frak F)}
\def\ZN{{\mathbb Z}_N}%
\def\1{\mathbbm 1}%
\def\module_1{{\frak C}_N^{\bar 1}}%
\def\A{{\mathscr{A}}}%
\def\N1{{\mathbb N}_1}%

\begin{document}

%
%
%
%
%
%
%
%
%

\title{Noncommutative Galois Extension and\\ Graded $q$-Differential Algebra}

\author{Viktor Abramov}

\address{%
Institute of Mathematics\\
University of Tartu, J. Liivi 2,\\
50409 Tartu, Estonia}

\email{viktor.abramov@ut.ee}

\subjclass{Primary 46L87; Secondary 81R60}

\keywords{Graded $q$-differential algebra, Galois noncommutative extension, reduced quantum plane}

\date{January 1, 2004}

\begin{abstract}
We show that a semi-commutative Galois extension of a unital associative algebra can be endowed with the structure of a graded $q$-differential algebra. We study the first and higher order noncommutative differential calculus of semi-commutative Galois extension induced by the graded $q$-differential algebra. As an example we consider the quaternions which can be viewed as the semi-commutative Galois extension of complex numbers. 
\end{abstract}

\maketitle
\section{Introduction}
Let $\tilde\A$ be an associative unital $\mathbb C$-algebra and $\A$ be its subalgebra. Assume that there is an element $\tau\in\tilde\A$ which does not belong to subalgebra $\A$ and $\tau^N=\1$, where $N\geq 2$ is an integer and $\1$ is the identity element of $\tilde\A$. A noncommutative Galois extension of $\A$ by means of $\tau$ is the smallest subalgebra $\A[\tau]\subset \tilde\A$ such that $\A\subset \A[\tau]$, and $\tau\in\A[\tau]$. A notion of noncommutative Galois extension of associative unital complex algebra was introduced and studied in the series of papers \cite{Kerner-Suzuki,Lawrynowicz-al-1,Lawrynowicz-al-2,Trovon}. It should be mentioned here that in the papers \cite{Lawrynowicz-al-2,Trovon} a concept of binary and ternary noncommutative Galois extension was studied in the relation with a ternary Clifford algebra, and it was shown there that these structures can be applied to theoretical elementary particle physics in order to construct an elegant algebraic model for quarks.

A concept of a graded $q$-differential algebra can be viewed as a result of development of idea of transition from the equation $d^2=0$ to the more general equation $d^N=0$, where $d$ is a differential of graded differential algebra and $N\geq 2$ is an integer. This idea was proposed and developed with the help of primitive $N$th root of unity in the paper \cite{Kapranov}, where the author introduced the notions of $N$-complex and generalized cohomologies of $N$-complex. Later this idea was developed in the paper \cite{Dubois-Violette-Kerner}, where the authors introduced and studied a notion of graded $q$-differential algebra. It was shown \cite{Abramov-1,Abramov-2,Abramov-Liivapuu,Abramov-Liivapuu-2} that a notion of a graded $q$-differential algebra can be applied in noncommutative geometry in order to construct a noncommutative generalization of differential forms and a concept of connection.

In this paper we apply the methods of graded $q$-differential algebras to a semi-commutative Galois extension of an associative unital $\mathbb C$-algebra. We show that a semi-commutative Galois extension can be endowed with the structure of graded $q$-differential algebra. We study the first and higher order noncommutative differential calculus induced by the $N$-differential of graded $q$-differential algebra. We introduce a derivative and differential with the help of first order noncommutative differential calculus developed in the papers \cite{Abramov-Kerner,Borowiec-Kharchenko}.





\section{Graded $q$-differential algebra structure of noncommutative Galois extension}
In this section our aim is to show that given a noncommutative Galois extension we can construct a graded $q$-differential algebra, where $q$ is a primitive $N$th root of unity.

\vskip.3cm
\noindent
Let $\A=\oplus_{k\in{\mathbb Z}_N}{\A}^k=\A^0\oplus\A^1\oplus\ldots\oplus\A^{N-1}$ be a ${\mathbb Z}_N$-graded associative unital $\mathbb C$-algebra with identity element denoted by $ \1$. Obviously the subspace $\A^0$ of elements of degree 0 is the subalgebra of a graded algebra $\A$.
Every subspace $\A^k$ of homogeneous elements of degree $k\geq 0$ can be viewed as the $\A^0$-bimodule.

\vskip.3cm
\noindent
Let us remind several basic notions of $q$-calculus, where $q$ is a complex number. The graded $q$-commutator of two homogeneous elements $u,v\in \A$ is defined  by
$$
[v,u]_q=v\,u-q^{|v||u|}u\,v.
$$
A graded $q$-derivation of degree $m$ of a graded algebra $\A$ is a linear mapping $d:\A\to\A$ of degree $m$, i.e. $d:\A^i\to \A^{i+m}$, which satisfies the graded $q$-Leibniz rule
\begin{equation}
d(u\,v)=d(u)\,v+q^{ml}u\,d(v),
\label{graded q-Leibniz}
\end{equation}
where $u$ is a homogeneous element of degree $l$, i.e. $u\in \A^l$. A graded $q$-derivation $d$ of degree $m$ is called an inner graded $q$-derivation of degree $m$ induced by an element $v\in\A^m$ if
\begin{equation}
d(u)=[v,u]_q=v\,u-q^{ml}u\,v,
\end{equation}
where $u\in \A^l$.

\vskip.3cm
\noindent
Now let $q$ be a primitive $N$th root of unity, for instant $q=e^{2\pi i/N}$. Then
$$
q^N=1,\;\; 1+q+\ldots+q^{N-1}=0.
$$
A graded $q$-differential algebra
is a graded associative unital algebra $\A$ endowed with a graded $q$-derivation $d$ of degree one which satisfies $d^N=0$. In what follows a graded $q$-derivation $d$ of a graded $q$-differential algebra $\A$ will be referred to as a graded $N$-differential. Thus a graded $N$-differential $d$ of a graded $q$-differential algebra is a linear mapping of degree one which satisfies a graded $q$-Leibniz rule and $d^N=0$. It is useful to remind that a graded differential algebra is a graded associative unital algebra equipped with a differential $d$ which satisfies the graded Leibniz rule and $d^2=0$. Hence it is easy to see that a graded differential algebra is a particular case of a graded $q$-differential algebra when $N=2, q=-1$, and in this sense we can consider a graded $q$-differential algebra as a generalization of a concept of graded differential algebra. Given a graded associative algebra $\A$ we can consider the vector space of inner graded $q$-derivations of degree one of this algebra and put the question: under what conditions an inner graded $q$-derivation of degree one is a graded $N$-differential? The following theorem gives  answer to this question.
\begin{thm}
Let $\A$ be a $\mathbb Z_N$-graded associative unital $\mathbb C$-algebra and $d(u)=[v,u]_q$ be its inner graded $q$-derivation induced by an element $v\in\A^1$. The inner graded $q$-derivation $d$ is the $N$-differential, i.e. it satisfies $d^N=0$, if and only if $v^N=\pm \1$.
\end{thm}
\noindent
Our aim in this section is to show that we can apply this theorem to a noncommutative Galois extension to construct a graded $q$-differential algebra with $N$-differential satisfying $d^N=0$. First of all we remind a notion of a noncommutative Galois extension \cite{Kerner-Suzuki,Lawrynowicz-al-1,Lawrynowicz-al-2,Trovon}.
\begin{defn}
Let $\tilde\A$ be an associative unital $\mathbb C$-algebra and $\A\subset \tilde\A$ be its subalgebra. If there exist an element $\tau\in\tilde\A$ and an integer $N\geq 2$ such that
\begin{itemize}
\item[i)] $\tau^N=\pm\1$,
\item[ii)] $\tau^k\notin\A$ for any integer $1\leq k\leq N-1$,
\end{itemize}
then the smallest subalgebra $\A[\tau]$ of $\tilde\A$ which satisfies
\begin{itemize}
\item[iii)] $\A\subset \A[\tau]$,
\item[iv)] $\tau\in \A[\tau]$,
\end{itemize}
is called the noncommutative Galois extension of $\A$ by means of $\tau$.
\end{defn}
\noindent
In this paper we will study a particular case of a noncommutative Galois extension which is called a semi-commutative Galois extension \cite{Trovon}. We will give a definition of a semi-commutative Galois extension with the help of left and right $\A$-modules generated by $\tau$. Let $\A^1_{\mbox{l}}[\tau]$ and $\A^1_{\mbox{r}}[\tau]$ be respectively the left and right $\A$-modules generated by $\tau$. Obviously we have
$$
\A^1_{\mbox{l}}[\tau]\subset \A[\tau],\;\;\A^1_{\mbox{r}}[\tau]\subset\A[\tau].
$$
\begin{defn}
A noncommutative Galois extension $\A[\tau]$ is said to be a right (left) semi-commutative Galois extension if $\A^1_{\mbox{r}}[\tau]\subset \A^1_{\mbox{l}}[\tau]$ ($\A^1_{\mbox{l}}[\tau]\subset \A^1_{\mbox{r}}[\tau]$). If $\A^1_{\mbox{r}}[\tau]\equiv \A^1_{\mbox{l}}[\tau]$ then a noncommutative Galois extension will be referred to as a semi-commutative Galois extension, and in this case $\A^1[\tau]=\A^1_{\mbox{r}}[\tau]=\A^1_{\mbox{l}}[\tau]$ is the $\A$-bimodule.
\end{defn}
\noindent
It is well known that a bimodule over an associative unital algebra $\A$ freely generated by elements of its basis induces the endomorphism from an algebra $\A$ to the algebra of square matrices over $\A$. In the case of semi-commutative Galois extension we have only one generator $\tau$ and it induces the endomorphism of an algebra $\A$. Indeed let $\A[\tau]$ be a semi-commutative Galois extension and $\A^1[\tau]$ be its $\A$-bimodule generated by $[\tau]$. Any element of the right $\A$-module $\A^1_{\mbox{r}}[\tau]$ can be written as $\tau\,x$, where $x\in \A$. On the other hand $\A[\tau]$ is a semi-commutative Galois extension which means $\A^1_{\mbox{r}}[\tau]\equiv \A^1_{\mbox{l}}[\tau]$, and hence each element $x\,\tau$ of the left $\A$-module can be expressed as $\tau\,\phi_\tau(x)$, where $\phi_\tau(x)\in\A$. It is easy to verify that the linear mapping $\phi:x\to \phi_\tau(x)$ is the endomorphism of subalgebra $\A$, i.e. for any elements $x,y\in \frak A$ we have $\phi_\tau(xy)=\phi_\tau(x)\phi_\tau(y)$. This endomorphism will play an important role in our differential calculus, and in what follows we will also use the notation $\phi_\tau(x)=x_\tau$. Thus
$$
u\,\tau=\tau\,\phi_\tau(x), \quad u\,\tau=\tau\,u_\tau.
$$
It is clear that
$$
\phi_\tau^N=\mbox{id}_{\A},\,\,\,u_{\tau^N}=u,
$$
because for any $u\in\A$ it holds $u\,\tau^N=\tau^N\,\phi^N(u)$ and taking into account that $\tau^N=\1$ we get $\phi_\tau^N(u)=u$.
\begin{prop}
Let $\A[\tau]$ be a semi-commutative Galois extension of $\A$ by means of $\tau$, and $\A_{\mbox{l}}^k[\tau],\A_{\mbox{r}}^k[\tau]$ be respectively the left and right $\A$-modules generated by $\tau^k$, where $k=1,2,\ldots, N-1$. Then  $\A_{\mbox{l}}^k[\tau]\equiv \A_{\mbox{r}}^k[\tau]=\A^k[\tau]$ is the $\A$-bimodule, and
$$
\A[\tau]=\oplus_{k=0}^{N-1}\A^k[\tau]=\A^0[\tau]\oplus\A^1[\tau]\oplus\ldots\oplus\A^{N-1}[\tau],
$$
where $\A^0[\tau]\equiv \A$.
\label{proposition Galois extension}
\end{prop}
\noindent
Evidently the endomorphism of $\A$ induced by the $\A$-bimodule structure of $A^k[\tau]$ is $\phi^k$, where $\phi:\A\to \A$ is the endomorphism induced by the $\A$-bimodule $\A^1[\tau]$. We will also use the notation $\phi^k(x)=x_{\tau^k}$.

\vskip.3cm
\noindent
It follows from Proposition \ref{proposition Galois extension} that a semi-commutative Galois extension $\A[\tau]$ has a natural $\ZN$-graded structure which can be defined as follows: we assign degree zero to each element of subalgebra $\A$, degree 1 to $\tau$ and extend this graded structure to a semi-commutative Galois extension $\A[\tau]$ by determining the degree of a product of two elements as the sum of degree of its factors. The degree of a homogeneous element of $\A[\tau]$ will be denoted by $|\;\;|$. Hence $|u|=0$ for any $u\in\A$ and $|\tau|=1$.
\begin{prop}
Let $q$ be a primitive $N$th root of unity. A semi-commutative Galois extension $\A[\tau]$, equipped with the $\ZN$-graded structure described above and with the inner graded $q$-derivation $d\,=[\tau,\,]_q$ induced by $\tau$, is the graded $q$-differential algebra, and $d$ is its $N$-differential. For any element $\xi$ of semi-commutative Galois extension $\A[\tau]$ written as a sum of elements of right $\A$-modules $\A^k[\tau]$
$$
\xi=\sum_{k=0}^{N-1}\tau^k\,u_k=\1\,u_0+\tau\,u_1+\tau^2\,u_2+\ldots \tau^{N-1}\,u_{N-1},\;\;u_k\in \A,
$$
it holds
\begin{equation}
d \xi=\sum_{k=0}^{N-1}\tau^{k+1}(u_k-q^k\,(u_k)_\tau),
\label{semi-commutative graded q-differential}
\end{equation}
where $u_k\to (u_k)_\tau$ is the endomorphism of $\A$ induced by the bimodule structure of $\A^1[\tau].$
\label{proposition 2.5}
\end{prop}
\section{First order differential calculus of noncommutative Galois extension}
In this section we develop a first order differential calculus of noncommutative Galois extension induced by the structure of graded $q$-differential algebra and its $N$-differential $d$. Let us remind that in noncommutative geometry a differential calculus is a triple $(\A,d,\mathscr M)$, where $\A$ is a unital associative algebra, $\mathscr M$ is an $\A$-bimodule, $d$ is a linear mapping $d:\A\to \mathscr M$ which satisfies Leibniz rule. If this calculus is a coordinate calculus \cite{Borowiec-Kharchenko} then one can introduce the partial derivatives induced by this calculus.

\vskip.3cm
\noindent
Let $\A[\tau]$ be a semi-commutative Galois extension of an algebra $\A$ by means of $\tau$. Thus we have an algebra $\A$ and $\A$-bimodule $\A^1[\tau]$. Next we have the $N$-differential $d:\A[\tau]\to\A[\tau]$ induced by $\tau$, and if we restrict this $N$-differential to the subalgebra $\A$ of Galois extension $\A[\tau]$ then $d:\A\to \A^1[\tau]$ satisfies the Leibniz rule. Consequently we have the first order differential calculus which can be written as the triple $(\A,d,\A^1[\tau])$. In order to describe the structure of this first order differential calculus we will need the vector space endomorphism $\Delta:\A\to\A$ defined by
$$
\Delta u=u-u_\tau,\quad u\in\A.
$$
For any elements $u,v\in\A$ this endomorphism satisfies
$$
\Delta(u\,v)=\Delta(u)\,v+u_\tau\,\Delta(v).
$$
Let us assume that there exists an element $x\in \A$ such that the element $\Delta x\in \A$ is invertible, and the inverse element will be denoted by $\Delta x^{-1}.$ The differential $dx$ of an element $x$ can be written in the form $dx=\tau\,\Delta x$ which clearly shows that $dx$ has degree one, i.e. $dx\in\A^1[\tau]$, and hence $dx$ can be used as generator for the right $\A$-module $\A^1[\tau]$. Let us denote by $\phi_{dx}:u\to \phi_{dx}(u)=u_{dx}$ the endomorphism of $\A$ induced by bimodule structure of $\A^1[\tau]$ in the basis $dx$. Then
\begin{equation}
u_{dx}=\Delta x^{-1}\,u_\tau\,\Delta x=\mbox{Ad}_{\Delta\,x}\,u_\tau.
\end{equation}
\begin{defn}
For any element $u\in \A$ we define the right derivative $\frac{du}{dx}\in \A$ (with respect to $x$)  by the formula
\begin{equation}
du=dx\,\frac{du}{dx}.
\label{derivative}
\end{equation}
\end{defn}
\noindent
Analogously one can define the left derivative with respect to $x$ by means of the left $\A$-module structure of $\A^1[\tau]$. Further we will only use the right derivative which will be referred to as the derivative and often will be denoted by $u^\prime_x$. Thus we have the linear mapping
$$
\frac{d}{dx}:\A\to\A, \quad \frac{d}{dx}:u\mapsto u^\prime_x.
$$
\begin{prop}
For any element $u\in \A$ we have
\begin{equation}
\frac{du}{dx}=\Delta x^{-1}\,\Delta u.
\end{equation}
The derivative (\ref{derivative}) satisfies the twisted Leibniz rule, i.e. for any two elements $u,v\in \A$ it holds
$$
\frac{d}{dx}(u\,v)=\frac{du}{dx}\,v+\phi_{dx}(u)\,\frac{dv}{dx}=\frac{du}{dx}\,v+\mbox{Ad}_{\Delta\,x}\,u_\tau\,\frac{dv}{dx}.
$$
\label{twisted Leibniz rule}
\end{prop}
\noindent
We have constructed the first order differential calculus with one variable $x$, and it is natural to study a transformation rule of the derivative of this calculus if we choose another variable. From the point of view of differential geometry we will study a change of coordinate in one dimensional space. Let $y\in\A$ be an element of $\A$ such that $\Delta\,y=y-y_\tau$ is invertible.
\begin{prop}
Let $x,y$ be elements of $\A$ such that $\Delta\,x,\Delta\,y$ are invertible elements of $\A$. Then
$$
dy=dx\,\,y^\prime_x,\;\; \frac{d}{dx}=y^\prime_x\,\frac{d}{dy},\;\;
     dx=dy\,\,x^\prime_y,\;\;\frac{d}{dy}=x^\prime_y\,\frac{d}{dx},
$$
where $x^\prime_y=(y^\prime_x)^{-1}$.
\end{prop}
\noindent
Indeed we have $dy=\tau\,\Delta\,y,\,dx=\tau\,\Delta\,x$. Hence $\tau=dx\,\Delta\,x^{-1}$ and
$$
dy=dx\,(\Delta\,x^{-1}\Delta\,y)=dx\,y^\prime_x.
$$
If $u$ is any element of $\A$ the for the derivatives we have
$$
\frac{du}{dx}=\Delta\,x^{-1}\,\Delta\,u=(\Delta\,x^{-1}\,\Delta\,y)\,(\Delta\,y^{-1}\,\Delta\,u)=y^\prime_x\,\frac{du}{dy}.
$$
As an example of the structure of graded $q$-differential algebra induced by $d_\tau$ on a semi-commutative Galois extension we can consider the quaternion algebra $\mathbb H$. The quaternion algebra $\mathbb H$ is associative unital algebra generated over $\mathbb R$ by $i,j,k$ which are subjected to the relations
$$
i^2=j^2=k^2=-\1,\;i\,j=-j\,i=k,\;j\,k=-k\,j=i,\;k\,i=-i\,k=j,
$$
where $\1$ is the unity element of $\H$. Given a quaternion
$$
\frak q=a_0\,\1+a_1\,i+a_2\,j+a_3\,k
$$
we can write it in the form $\frak q=(a_0\,\1+a_2\,j)+i\,(a_1+a_3\,j)$. Hence if we consider the coefficients of the previous expression $z_0=a_0\,\1+a_2\,j,z_1=a_1+a_3\,j$ as complex numbers then $\frak q=z_0\,\1+i\,z_1$ which clearly shows that the quaternion algebra $\H$ can be viewed as the semi-commutative Galois extension $\mathbb C[i]$. Evidently in this case we have $N=2,q=-1$, and $\mathbb Z_2$-graded structure defined by $|\1|=0, |i|=1$. Hence we can use the terminology of superalgebras. It is easy to see that the subspace of odd elements (degree 1) can be considered as the bimodule over the subalgebra of even elements $a\,\1+b\,j$ and this bimodule induces the endomorphism $\phi:\mathbb C\to\mathbb C$, where $\phi(z)=\bar{z}$. Let $d$ be the differential of degree one (odd degree operator) induced by $i$. Then making use of \ref{semi-commutative graded q-differential} for any quaternion $\frak q$ we have
$$
d\frak q=d(z_0\,\1+i\,z_1)=-(\bar{z}_1+z_1)\,\1.
$$
Obviously $d^2\frak q=0.$
\section{Higher order differential calculus}
Our aim in this section is to develop a higher order differential calculus of a semi-commutative Galois extension $\A[\tau]$. This higher order differential calculus is induced by the graded $q$-differential algebra structure. In Section 2 it is mentioned that a graded $q$-differential algebra can be viewed as a generalization of a concept of graded differential algebra if we take $N=2,q=-1$. It is well known that one of the most important realizations of graded differential algebra is the algebra of differential forms on a smooth manifold. Hence we can consider the elements of the graded $q$-differential algebra constructed by means of a semi-commutative Galois extension $\A[\tau]$ and expressed in terms of differential $dx$ as noncommutative analogs of differential forms with exterior differential $d$ which satisfies $d^N=0$. In order to stress this analogy we will consider an element $x\in\A$ as analog of coordinate, the elements of degree zero as analogs of functions, elements of degree $k$ as analogs of $k$-forms, and we will use the corresponding terminology. It should be pointed out that because of the equation $d^N=0$ there are higher order differentials $dx,d^2x,\ldots,d^{N-1}x$ in this algebra of differential forms.

Before we describe the structure of higher order differentials forms it is useful to introduce the polynomials $P_k(x),Q_k(x)$, where $k=1,2,\ldots,N$. Let us remind that $\Delta x=x-x_\tau\in\A$. Applying the endomorphism $\tau$ we can generate the sequence of elements
$$
\Delta x_\tau=x_\tau-x_{\tau^2}, \Delta x_{\tau^2}=x_{\tau^2}-x_{\tau^3},\ldots, \Delta x_{\tau^{N-1}}=x_{\tau^{N-1}}-x.
$$
Obviously each element of this sequence is invertible. Now we define the sequence of polynomials $Q_1(x),Q_2(x),\ldots,Q_N(x)$, where
$$
Q_k(x)=\Delta x_{\tau^{k-1}}\Delta x_{\tau^{k-2}}\ldots \Delta x_{\tau}\Delta x.
$$
These polynomials can be defined by means of the recurrent relation
$$
Q_{k+1}(x)=(Q_k(x))_\tau\Delta x.
$$
It should be mentioned that $Q_k(x)$ is the invertible element and
$$
(Q_k(x))^{-1}=\Delta x^{-1}\Delta x^{-1}_\tau\ldots \Delta x^{-1}_{\tau^{k-1}}.
$$
We define the sequence of elements $P_1(x),P_2(x),\ldots,P_N(x)\in\A$ by the recurrent formula
$$
P_{k+1}(x)=P_{k}(x)-q^{k}\,(P_{k}(x))_\tau,\quad k=1,2,\ldots,N-1,
$$
and $P_1(x)=\Delta x$. Clearly $P_1(x)=Q_(x)$ and for the $k=2,3$ a straightforward calculation gives
\begin{eqnarray}
P_2(x) &=& x-(1+q)\,x_\tau+q\,x_{\tau^2},\nonumber\\
P_3(x) &=& x-(1+q+q^2)\,x_\tau+(q+q^2+q^3)\,x_{\tau^2}-q^3\,x_{\tau^3}.\nonumber
\end{eqnarray}
\begin{prop}
If $q$ is a primitive $N$th root of unity then there are the identities
$$
P_{N-1}(x)+(P_{N-1}(x))_\tau+\ldots+(P_{N-1}(x))_{\tau^{N-1}}\equiv 0,\;\;\; P_N(x)\equiv 0.
$$
\end{prop}
\noindent
Now we will describe the structure of higher order differential forms. It follows from the previous section that any 1-form $\omega$, i.e. an element of $\A^1[\tau]$, can be written in the form $\omega=dx\,u$, where $u\in\A$. Evidently $d:\A\to \A^1[\tau],\,d\omega=dx\,u^{\prime}_x$. The elements of $\A^2[\tau]$ will be referred to as 2-forms. In this case there are two choices for a basis for the right $\A$-module $\A^2[\tau]$. We can take either $\tau^2$ or $(dx)^2$ as a basis for $\A^2[\tau]$. Indeed we have
$$
(dx)^2=\tau^2\,Q_2(x).
$$
It is worth mentioning that the second order differential $d^2x$ can be used as the basis for $\A^2[\tau]$ only in the case when $P_2(x)$ is invertible. Indeed we have
$$
d^2x=\tau^2\,P_2(x),\quad d^2x=(dx)^2\,Q^{-1}_2(x)P_2(x).
$$
If we choose $(dx)^2$ as the basis for the module of 2-forms $\A^2[\tau]$ then any 2-form $\omega$ can be written as $\omega=(dx)^2\,u$, where $u\in\A$. Now the differential of any 1-form $\omega=dx\,u$, where $u\in\A$, can be expressed as follows
\begin{equation}
d\omega=(dx)^2\,\big( q\,u^\prime_x+Q_2^{-1}(x)P_2(x)\,u \big).
\label{differential of 1-form}
\end{equation}
It should be pointed out that the second factor of the right-hand side of the above formula resembles a covariant derivative in classical differential geometry. Hence we can introduce the linear operator $D:\A\to \A$ by the formula
\begin{equation}
Du=q\,u^\prime_x+Q_2^{-1}(x)P_2(x)\,u,\quad u\in\A.
\label{covariant derivative 1}
\end{equation}
If $\omega=dv, v\in \A$, i.e. $\omega$ is an exact form, then
$$
d\omega=d^2v=(dx)^2\,Dv^\prime_x=(dx)^2\,\big(q\,v^{\prime\prime}_x+Q_2^{-1}(x)P_2(x)\,v_x^\prime\big).
$$
If we consider the simplest case $N=2,q=-1$ then
$$
d^2v=0,\;\; P_2(x)\equiv 0,\;\; (dx)^2\neq 0,
$$
and from the above formula it follows that $v^{\prime\prime}_x=0$.
\begin{prop}
Let $\A[\tau]$ be a semi-commutative Galois extension of algebra $\A$ by means of $\tau$, which satisfies $\tau^2=\1$, and $d$ be the differential of the graded differential algebra induced by an element $\tau$ as it is shown in Proposition \ref{proposition 2.5}. Let $x\in\A$ be an element such that $\Delta x$ is invertible. Then for any element $u\in\A$ it holds $u^{\prime\prime}_x=0$, where $u^\prime_x$ is the derivative (\ref{derivative}) induced by $d$. Hence any element of an algebra $\A$ is  linear with respect to $x$.
\label{linear}
\end{prop}
\noindent
The quaternions considered as the noncommutative Galois extension of complex numbers (Section 3) provides a simple example for the above proposition. Indeed in this case $\tau=i, \A\equiv \mathbb C$, where the imaginary unit is identified with $j$, $(a\,\1+b\,j)_\tau=a\,\1-b\,j$. Hence we can choose $x=a\,\1+b\,j$ iff $b\neq 0$. Indeed in this case $\Delta x=x-x_\tau=a\,\1+b\,j-a\,\1+b\,j=2b\,j$, and $\Delta x$ is invertible iff $b\neq 0$. Now any $z=c\,\1+d\,j\in\A$ can be uniquely written in the form $z=\tilde c\,\1+\tilde d\,x$ iff
$$
\left|
  \begin{array}{cc}
    1 & a \\
    0 & b \\
  \end{array}
\right|=b\neq 0.
$$
Thus any $z\in\A$ is linear with respect to $x$.

Now we will describe the structure of module of $k$-forms $\A^k[\tau]$. We choose $(dx)^k$ as the basis for the right $\A$-module $\A^k[\tau]$, then any $k$-form $\omega$ can be written $\omega=(dx)^k\,u, \;u\in\A$. We have the following relations
$$
(dx)^k=\tau^k\,Q_k(x),\;\;\;d^kx=\tau^k\,P_k(x).
$$
In order to get a formula for the exterior differential of a $k$-form $\omega$ we need the polynomials $\Phi_1(x),\Phi_2(x),\ldots,\Phi_{N-1}(x)$ which can be defined by the recurrent relation
\begin{equation}
\Phi_{k+1}(x)=\text{Ad}_{\Delta x}(\Phi_{k}) + q^{k-1}\Phi_1(x),\quad k=1,2,\ldots,N-1,
\end{equation}
where $\Phi_1(x)=Q^{-1}_2(x)P_2(x)$. These polynomials satisfy the relations $d(dx)^k=(dx)^{k+1}\Phi_k(x)$ and given a $k$-form $\omega=(dx)^k\,u, \;u\in\A$ we find its exterior differential as
$$
d\omega=(dx)^{k+1}\,\bigg(q^k\,u^\prime_x + \Phi_k(x)\,u\bigg)=(dx)^{k+1}\,D^{(k)}u.
$$
The linear operator $D^{(k)}:\A\to\A, k=1,2,\ldots,N-1$ introduced in the previous formula has the form
\begin{equation}
D^{(k)}u=q^k\,u^\prime_x + \Phi_k(x)\,u,
\label{covariant derivative 2}
\end{equation}
and, as it was mentioned before, this operator resembles a covariant derivative of classical differential geometry. It is easy to see that the operator (\ref{covariant derivative 1}) is the particular case of (\ref{covariant derivative 2}), i.e. $D^{(1)}\equiv D$.
\section{Conclusions}
A graded $q$-differential algebra structure arises in a very natural way in the case of a semi-commutative Galois extension. The first order differential calculus of this algebra allows us to introduce the derivative and construct the noncommutative analog of differential forms with a differential satisfying $d^N=0$. Though we have an analogy with differential forms it should be mentioned that the calculus of differential forms with a differential $d^N=0$ considered in the present paper has some differences from the classical differential forms on a smooth manifold. First of all if we consider the simplest case of our calculus when $N=2,q=-1$ then $d^2=0$ not because of skew-commutativity of differentials in a wedge product (and the property of second order derivatives which do not depend on the order of differentiation), but simply because of the fact that the second order derivative of any function of our calculus is zero, i.e. we have only linear functions (Proposition \ref{linear}). Secondly a peculiar property of our calculus is the appearance of higher order differentials of "coordinate" $dx,d^2x,\ldots,d^{N-1}x$. As it was mentioned before in the case of these higher order differentials a differential $d$ induces the operator (\ref{covariant derivative 1},\ref{covariant derivative 2}) which is very similar to a covariant derivative of classical differential geometry. It is worth to mention here that an analogous result was obtained in the paper \cite{Abramov-Kerner}, where the authors constructed analogs of differential forms on a smooth manifold with exterior differential $d^3=0$, and it was shown that in the case of second order differentials one should use a covariant derivative.

A reduced quantum plane can be viewed as the algebra generated by $x,y$ which obey the relations $x\,y=q\;y\,x, \; x^N=y^N=\1$. It should be pointed out that this algebra is the particular case of a generalized Clifford algebra. Obviously a reduced quantum plane can be considered as the semi-commutative Galois extension of the algebra generated by $x$ by means of $y$. We intend to apply the approach developed in the present paper to a reduced quantum plane to study it with the help of graded $q$-differential algebra and the corresponding calculus of analogs of differential forms.

\subsection*{Acknowledgment}
The authors is gratefully acknowledge the Estonian Science Foundation for financial support of this work under the Research Grant No. ETF9328. This research was also supported by institutional research funding IUT20-57
of the Estonian Ministry of Education and Research.

\end{document}